\documentclass[oneside,english]{amsart}
\usepackage[T1]{fontenc}
\usepackage[latin1]{inputenc}
\usepackage{amssymb}

\makeatletter
\theoremstyle{plain}
\newtheorem{thm}{Theorem}[section]
\theoremstyle{plain}
\newtheorem{prop}[thm]{Proposition}
\theoremstyle{plain}
\newtheorem{cor}[thm]{Corollary}
\theoremstyle{remark}
\newtheorem*{rem*}{Remark}
\theoremstyle{remark}
\newtheorem*{acknowledgement*}{Acknowledgement}

\usepackage{babel}
\makeatother
\begin{document}

\title[Geometry of systems of third order
differential equations]{The geometry of systems of third order
differential equations induced by second order Lagrangians}
\author{Ioan Bucataru and Radu Miron}
\address{Faculty of Mathematics, "Al.I.Cuza" University, Iasi,
700506, Romania} \email{bucataru@uaic.ro, radu.miron@uaic.ro}
\urladdr{http://www.math.uaic.ro/~bucataru}

\begin{abstract}
A dynamical system on the total space of the fibre bundle of second
order accelerations, $T^2M$, is defined as a third order vector
field $S$ on $T^2M$, called semispray, which is mapped by the second
order tangent structure into one of the Liouville vector field. For
a regular Lagrangian of second order we prove that this semispray is
uniquely determined by two associated Cartan-Poincar\'e one-forms.
To study the geometry of this semispray we construct a nonlinear
connection, which is a Lagrangian subbundle for the presymplectic
structure. Using this semispray and the associated nonlinear
connection we define covariant derivatives of first and second
order. With respect to this, the second order dynamical derivative
of the Lagrangian metric tensor vanishes.
\end{abstract}

\maketitle



\section{Introduction}

The notion of spray as a vector field that lives on the total
space of the first order tangent bundle was introduced by Ambrose
et al. \cite{ambrose60}. For a regular Lagrangian of order one,
solutions of the Euler-Lagrange equations are integral curves of
the canonical semispray, which is a second order vector field. The
canonical semispray of a Lagrange space is uniquely determined by
its symplectic structure and the energy of the space. This point
of view allows for studying first integrals, Cartan symmetries,
and Noether type theorems and it has been investigated by Klein
\cite{klein62}, Godbillon \cite{godbillon69}, Abraham and Marsden
\cite{abraham78}, Arnold \cite{arnold78}, Crampin and Pirani
\cite{crampin_pirani86}, de León and Rodrigues \cite{deLeon89},
Krupkova \cite{krupkova1}, Miron and Anastasiei \cite{mirona94}.
For the differential geometry of a second order vector field one
has to study its associated nonlinear connection and dynamical
covariant derivative. Such nonlinear connection was introduced by
Crampin \cite{crampin71} and Grifone \cite{grifone72} and it is a
very important tool for the geometry of Finsler and Lagrange
spaces developed by Miron and Anastasiei \cite{mirona94}. For a
regular Lagrangian, the canonical nonlinear connection is the only
connection that is metric and compatible with the symplectic
structure, as it has been shown in \cite{bucataru05}.

Tangent bundles of higher order are canonically endowed with tangent
structures of higher order, introduced by Eliopoulos
\cite{eliopoulos65}. This tangent geometry of higher order provides
a special differential calculus, which has been investigated by many
authors, for example by Yano and Ishihara \cite{yano76}, Tulczyjew
\cite{tulczyjew76}, de León and Rodrigues \cite{deLeon85}, Saunders
\cite{saunders89}, Miron \cite{miron97}. Important works that
connect the geometry of higher-order differential equations with the
geometry of higher order Lagrangian mechanics are due to Crampin et
al. \cite{cantrijn86}, de León and Rodrigues \cite{deLeon92}, Byrnes
\cite{byrnes96}.

Spaces involving metric tensors whose components are functions of
position and higher order accelerations were introduced first by
Kawaguchi \cite{kawaguchi32}. The variational problem of such metric
structures of order $k$, $k\geq 1$, was investigated by Craig
\cite{craig35} and Synge \cite{synge35} who associated to it a
number of $k+1$ covariant vectors. Among these vectors, only the
last two live on the tangent space of order $k$. For a regular
Lagrangian of order $k$, the last one determines one of the
Cartan-Poincar\'e one forms, the one before the last one determines
a vector field of order $k+1$. This vector field was called by Miron
in \cite{miron94} the Craig-Synge vector field, or the canonical
semispray, of the Lagrange space of order $k$. The geometry of
higher order Lagrangians based on the associated semispray and the
associated nonlinear connection has been developed recently by Miron
\cite{miron97}. From a semispray that lives on the tangent bundle of
order $k$, one can derive different nonlinear connections as it has
been done by Catz \cite{catz74}, de León and Rodrigues
\cite{deLeon85}, Byrnes \cite{byrnes96}, Miron and Atanasiu
\cite{mirona94}, Bucataru \cite{bucataru98}.  For $k>1$ these
nonlinear connections are different and the expression of their
coefficients is quite complicated. Geometric invariants induced by
two different nonlinear connections associated to such a semispray
were studied by Crampin and Saunders \cite{crampin06}. In this paper
we limit our work to the case of second order Lagrangians.

A dynamical system on the total space of second order acceleration
bundle $T^2M$ is defined as third order vector field, $S$, which
lives on $T^2M$ and is mapped by the second order tangent map into
one of the Liouville vector field. Such vector field is called a
semispray. For a regular Lagrangian of second order, we prove that
one can associate to it two Cartan-Poincar\'e one-forms and two
Cartan-Poincar\'e two-forms. We show that one of the
Cartan-Poincar\'e two-forms is a presymplectic form. Moreover, both
vertical distributions are Lagrangian subbundles with respect to
this presymplectic form.

The classical variational problem for a second order Lagrangian
leads to fourth-order equations, which are the Euler-Lagrange
equations. To the fourth-order Euler-Lagrange vector field, Craig
\cite{craig35} and Synge \cite{synge35} associated another vector
field of third order, which was called the Craig-Synge vector by
Miron \cite{miron94}. In this paper we study Craig-Synge equations
and the corresponding vector field since they live on the second
order tangent bundle of the base manifold. For this vector field,
which is a semispray of the considered regular Lagrangian of second
order, we prove that it is uniquely determined by the two associated
Cartan-Poincar\'e one-forms. This result extends the work of Sarlet
et al. \cite{sarlet84} for the canonical second order vector field
of a first order regular Lagrangian. We prove that the Craig-Synge
covector can be obtained directly by a special variational principle
of a second order regular Lagrangian for the case when only the
vertical part of the curve is varied. It is known that for a
Lagrangian of order grater than one, the integral of action along a
curve does not depend on the parametrization of the curve only if
the Lagrangian is degenerate. This is known as Zermelo condition,
geometric aspects of this condition were investigated by Kondo,
\cite{kondo63}. Therefore for a regular Lagrangian, when studying
the integral of action, the parametrization of a curve is essential.
The variation of the vertical part of the curve allows for variation
of such parametrization. Recently, seismic ray path through higher
order metric spaces, where the Zermelo condition is satisfied, were
investigated by Yajima and Nagahama \cite{yajima07}.

A nonlinear connection on the second order tangent bundle together
with a semispray define dynamical covariant derivatives of first and
second order. These covariant derivatives allow for invariant
expressions of various geometric objects derived from a second order
Lagrangian. In this paper we give conditions that uniquely determine
a new nonlinear connection of a second order regular Lagrangian. One
of these conditions states that the corresponding horizontal
distribution is a Lagrangian subbundle for the presymplectic form.
With respect to this nonlinear connection, the second order
dynamical covariant derivative of the metric tensor vanishes. This
generalizes the result for a first order regular Lagrangian, where
the dynamical covariant derivative of the metric tensor vanishes,
\cite{bucataru05}.

\section{Second order tangent bundle}

The space of second order accelerations can be identified with the
space of second order tangent bundle. The space of higher order
tangent bundles using jet theory was introduced by Ehresmann
\cite{ehresmann51}. In this section we present some of the geometric
structures that live on the total space of the second order tangent
bundle, such as second order tangent structures, Liouville vector
fields, vertical and horizontal distributions and semispray. These
structures were investigated previously by Eliopoulos
\cite{eliopoulos65}, Yano and Ishihara \cite{yano76}, Tulczyjew
\cite{tulczyjew76}, de León and Rodrigues \cite{deLeon85}, Crampin
et al. \cite{cantrijn86} and Miron \cite{miron97}.

For a real, $n$-dimensional, smooth manifold $M$, we denote by
$\left(T^{2}M,\pi^{2},M\right)$ its tangent bundle of second order.
A local chart $\left(U,\varphi=\left(x^{i}\right)\right)$ on $M$
induces a local chart
$\left(\left(\pi^{2}\right)^{-1}\left(U\right),
\Phi=\left(x^{i},y^{\left(1\right)i},y^{\left(2\right)i}\right)\right)$
on $T^{2}M$, where for a two-jet $j_{0}^{2}\rho\in T^{2}M$, the
coordinate functions are defined as follows \[
\left.y^{\left(\alpha\right)i}\left(j_{0}^{2}\rho\right)=
\frac{1}{\alpha!}\frac{d^{\alpha}\left(x^{i}\circ\rho\right)}{dt^{\alpha}}\right|_{t=0},
\  \alpha\in\left\{ 1,2\right\} .\] The second order tangent
structure $J$ and second order cotangent structure $J^{*}$, which
are called also vertical endomorphisms, are defined as follows\[
J=\frac{\partial}{\partial y^{\left(1\right)i}}\otimes
dx^{i}+\frac{\partial}{\partial y^{\left(2\right)i}}\otimes
dy^{\left(1\right)i},\ \ J^{*}=dx^{i}\otimes\frac{\partial}{\partial
y^{\left(1\right)i}}+dy^{\left(1\right)i}\otimes\frac{\partial}{\partial
y^{\left(2\right)i}}.\] The foliated structure of $T^{2}M$ allows
for two regular vertical distributions, $V_{1}\left(u\right)=Ker\
d_{u}\pi^{2}=Ker\ J^{2}=Im\ J$ and
$V_{2}\left(u\right)=J\left(V_{1}\left(u\right)\right)=Ker\
d_{u}\pi_{1}^{2}=Ker\ J=Im\ J^{2},$ $\forall u\in T^{2}M$. These two
distributions are integrable, the first one is tangent to the fibers
of $\pi^{2}:\left(x,y^{\left(1\right)},y^{\left(2\right)}\right)\in
T^{2}M\mapsto\left(x\right)\in M$, while the second one is tangent
to the fibers of
$\pi_{1}^{2}:\left(x,y^{\left(1\right)},y^{\left(2\right)}\right)\in
T^{2}M\mapsto\left(x,y^{\left(1\right)}\right)\in TM$.

The following two vertical vector fields are globally defined on $T^{2}M$
and they are called Liouville vector fields.

\[
{\mathbb C}_{2}= y^{\left(1\right)i}\frac{\partial}{\partial
y^{\left(1\right)i}} +2y^{\left(2\right)i}\frac{\partial}{\partial
y^{\left(2\right)i}},\ {\mathbb C}_{1}=J\left({\mathbb
C}_{2}\right)=y^{\left(1\right)i}\frac{\partial}{\partial
y^{\left(2\right)i}}.\] A semispray is a globally defined vector
field $S$ on $T^{2}M$ that satisfies the equation $JS={\mathbb
C}_{2}$. Therefore, a semispray $S$, which is a third order vector
field, can be expressed as follows

\begin{equation}
S=y^{\left(1\right)i}\frac{\partial}{\partial x^{i}}
+2y^{\left(2\right)i}\frac{\partial}{\partial
y^{\left(1\right)i}}-3G^{i}\left(x,y^{\left(1\right)},
y^{\left(2\right)}\right)\frac{\partial}{\partial
y^{\left(2\right)i}}\label{eq:sgi}\end{equation} and it is perfectly
determined by its coefficient functions
$G^{i}\left(x,y^{\left(1\right)},y^{\left(2\right)}\right)$. In this
paper we provide conditions that uniquely determine such a semispray
for a regular second order Lagrange space.

A nonlinear connection, or a horizontal distribution, on $T^{2}M$ is
a regular distribution $N: u\in T^{2}M\mapsto N\left(u\right)\subset
T_{u}T^{2}M$ that is supplementary to the vertical distribution
$V_{1}$. In other words, the following direct sum holds true.
\begin{equation} T_{u}T^{2}M=N\left(u\right)\oplus
V_{1}\left(u\right),\forall u\in
T^{2}M.\label{eq:dsum}\end{equation}

Since $d_{u}\pi^{2}:T_{u}T^{2}M\mapsto T_{\pi^{2}\left(u\right)}M$
is an epimorphism of vector spaces, whose kernel is
$V_{1}\left(u\right)$, its restriction to $N\left(u\right)$ is an
isomorphism of vector spaces. We denote by
$l_{h,u}:T_{\pi^{2}\left(u\right)}M\longrightarrow N\left(u\right)$
its inverse, to which we refer to as the horizontal lift induced by
the nonlinear connection $N$. If we consider the regular,
n-dimensional distribution
$N_{1}\left(u\right)=J\left(N\left(u\right)\right)$ then we have the
following direct sums. \begin{equation}
T_{u}T^{2}M=N\left(u\right)\oplus N_{1}\left(u\right) \oplus
V_{2}\left(u\right), \ V_{1}\left(u\right)=N_{1}\left(u\right)
\oplus V_{2}\left(u\right),\ \forall u\in
T^{2}M.\label{eq:dsum2}\end{equation} An adapted basis to the first
decomposition \eqref{eq:dsum2} is given by \begin{equation}
\left.\frac{\delta}{\delta
x^{i}}\right|_{u}:=l_{h,u}\left(\left.\frac{\partial}{\partial
x^{i}}\right|_{\pi^{2}\left(u\right)}\right),\left.\frac{\delta}{\delta
y^{\left(1\right)i}}\right|_{u}:=J\left(\left.\frac{\delta}{\delta
x^{i}}\right|_{u}\right),\left.\frac{\partial}{\partial
y^{\left(2\right)i}}\right|_{u}.\label{eq:adapted}\end{equation}
With respect to the natural basis of $T_{u}T^{2}M$, we use the
following notations:\begin{equation} \left.\frac{\delta}{\delta
x^{i}}\right|_{u}=\left. \frac{\partial}{\partial
x^{i}}\right|_{u}-\underset{\left(1\right)}{N_{i}^{j}}\left(u\right)\left.
\frac{\partial}{\partial y^{\left(1\right)j}}\right|_{u} -
\underset{\left(2\right)}{N_{i}^{j}}\left(u\right)\left.
\frac{\partial}{\partial y^{\left(2\right)j}}\right|_{u}.
\label{eq:delta}\end{equation} Functions
$\underset{\left(1\right)}{N_{i}^{j}}$ and
$\underset{\left(2\right)}{N_{i}^{j}}$ are called local coefficients
of the nonlinear connection $N$. The dual basis of the adapted basis
given by expression \eqref{eq:adapted} is given by \begin{equation}
dx^{i},\ \delta y^{\left(1\right)i}
=dy^{\left(1\right)i}+\underset{\left(1\right)}{M_{j}^{i}\ }dx^{j},\
\delta y^{\left(2\right)i}
=dy^{\left(2\right)i}+\underset{\left(1\right)}{M_{j}^{i}\ }
dy^{\left(1\right)j}+\underset{\left(2\right)}{M_{j}^{i}\ }dx^{j}.
\label{eq:dualbasis}\end{equation}

Functions $\underset{\left(1\right)}{M_{j}^{i}\ }$ and
$\underset{\left(2\right)}{M_{j}^{i}\ }$ are called the dual
coefficients of the nonlinear connection and they are related to the
coefficients of the nonlinear connection through the following
formulas\begin{equation} \underset{\left(1\right)}{M_{j}^{i}\
}=\underset{\left(1\right)}{N_{j}^{i}\ }, \
\underset{\left(2\right)}{M_{j}^{i}\
}=\underset{\left(2\right)}{N_{j}^{i}\
}+\underset{\left(1\right)}{N_{m}^{i}\
}\underset{\left(1\right)}{N_{j}^{m}\
}.\label{eq:nmrelation}\end{equation}

A semispray induces a nonlinear connection, which can be expressed
in terms of some product structures as it has been shown by Catz
\cite{catz74} and de León and Rodrigues \cite{deLeon85}. Different
sets of dual coefficients of nonlinear connections induced by a
semispray were obtained by Miron \cite{miron97} and Bucataru
\cite{bucataru98}. In this paper we provide conditions that uniquely
determine a nonlinear connection in terms of its compatibility with
the metric structure and the presymplectic structure of a second
order Lagrange space.

\section{Second order Lagrange space}

Important geometric structures for the geometry of a regular first
order Lagrangian can be derived from the associated
Cartan-Poincar\'e one and two-forms: Euler-Lagrange vector field,
nonlinear connection, symmetries. In this section we introduce some
of the geometric structures one can associate to a second order
regular Lagrangian, such as Cartan-Poincar\'e one and two-forms. We
study their relations with some other geometric structures that live
on the second order tangent bundle.

Consider $L:T^{2}M\longrightarrow\mathbb{R}$ a regular Lagrangian of
second order. In other words, the metric tensor \begin{equation}
g_{ij}=\frac{1}{2}\frac{\partial^{2}L}{\partial
y^{\left(2\right)i}\partial y^{\left(2\right)j}}
\label{eq:gij}\end{equation} is a symmetric second order tensor
field that has maximal rank $n$ on $T^{2}M.$

A tensor field on $T^{2}M$ is called a $d$-tensor field if its components
change as the components of a similar tensor field on the base manifold,
under corresponding change of coordinates. The metric tensor given
by expression \eqref{eq:gij} is a second order $d$-tensor field
since its components $g_{ij}$ behave as the components of a $\left(0,2\right)$-type
tensor field on the base manifold.

For a regular Lagrangian of second order, one can define the
following globally defined two Cartan-Poincar\'e one-forms:
\begin{equation} \theta_{L}^{1}=J^{*}\left(dL\right)=d_{J}L=
\frac{\partial L}{\partial y^{\left(1\right)i}}dx^{i} +
\frac{\partial L}{\partial y^{\left(2\right)i}}dy^{\left(1\right)i}
\ \textrm{\ and \ } \label{eq:thetal1}\end{equation}

\begin{equation}
\theta_{L}^{2}=\left(J^{*}\right)^{2}\left(dL\right)=
d_{J^{2}}L=\frac{\partial L}{\partial y^{\left(2\right)i}}dx^{i}.
\label{eq:thetal2}\end{equation} We consider also the following
Cartan-Poincar\'e two-forms: \begin{eqnarray} \omega_{L}^{2} & = &
d\theta_{L}^{2}= d\left(\frac{\partial L}{\partial
y^{\left(2\right)i}}dx^{i}\right)
= d\left(\frac{\partial L}{\partial y^{\left(2\right)i}}\right) \wedge dx^{i}=\nonumber \\
&  & \frac{\partial^{2}L}{\partial x^{j}\partial
y^{\left(2\right)i}}dx^{j}\wedge
dx^{i}+\frac{\partial^{2}L}{\partial y^{\left(1\right)j}\partial
y^{\left(2\right)i}}dy^{\left(1\right)j}\wedge
dx^{i}+2g_{ji}dy^{\left(2\right)j}\wedge
dx^{i},\label{eq:omegal2}\end{eqnarray}
\begin{equation}
\omega_{L}^{1}=d\theta_{L}^{1} = d\left(\frac{\partial L}{\partial
y^{\left(1\right)i}}\right)\wedge dx^{i}+d\left(\frac{\partial
L}{\partial y^{\left(2\right)i}}\right)\wedge
dy^{\left(1\right)i}.\label{eq:omegal1}\end{equation} We remark here
that the regularity of the Lagrangian $L$ implies the fact that
$rank\left(\omega_{L}^{2}\right)=2n<3n=\dim\left(T^{2}M\right)$. We
refer to $\omega_{L}^{2}$ as to the canonical presymplectic
structure of the Lagrangian $L$. These aspects were briefly
discussed in \cite{bucataru07}.

\begin{prop}
\label{pro:thetap} Consider $S$ an arbitrary semispray on $T^{2}M$.
For a second order Lagrangian $L$, its Cartan-Poincar\'e one and
two-forms have the following properties.
\end{prop}
\begin{enumerate}
\item $d_{J}\theta_{L}^{2}=d_{J}d_{J^{2}}L=0,$
$d_{J^{2}}\theta_{L}^{1}=d_{J^{2}}d_{J}L=0,$\medskip{}

\item $i_{S}\theta_{L}^{1}=\theta_{L}^{1}\left(S\right)={\mathbb C}_{2}\left(L\right),$
$i_{S}\theta_{L}^{2}=\theta_{L}^{2}\left(S\right)={\mathbb
C}_{1}\left(L\right),$\medskip{}

\item $\mathcal{L}_{S}\theta_{L}^{1}=i_{S}\omega_{L}^{1}
+d\left({\mathbb C}_{2}\left(L\right)\right),$
$\mathcal{L}_{S}\theta_{L}^{2}=i_{S}\omega_{L}^{2}+d\left({\mathbb
C}_{1}\left(L\right)\right).$
\end{enumerate}
\begin{proof}
First two properties follow directly from expressions
\eqref{eq:omegal2} and \eqref{eq:omegal1}, if we compose with $J$
and $J^{2}$ respectively. Next two properties are direct
consequences of expressions \eqref{eq:thetal1} and
\eqref{eq:thetal2}. Last two properties follow from the previous
ones if we take the exterior differential and use the fact that
$\mathcal{L}_{S}=di_{S}+i_{S}d,$ where $\mathcal{L}_{S}$ is the Lie
derivative along the semispray $S$.
\end{proof}
First property of Proposition \ref{pro:thetap} means that
$\omega_{L}^{2}\left(JX,JY\right)
=\omega_{L}^{2}\left(J^{2}X,J^{2}Y\right)=0,$ $\forall
X,Y\in\chi\left(T^{2}M\right)$ and therefore both vertical
distributions $V_{1}$ and $V_{2}$ are Lagrangian subbundles for the
presymplectic structure $\omega_{L}^{2}$. Second property of
Proposition \ref{pro:thetap} means that
$\omega_{L}^{1}\left(J^{2}X,J^{2}Y\right)=0,$ $\forall
X,Y\in\chi\left(T^{2}M\right)$ and therefore the vertical
distribution $V_{2}$ is a Lagrangian subbundle for the
Cartan-Poincar\'e two-form $\omega_{L}^{1}$.

We will explain now how the results in this section extend the
well-known results for the firs order case, to which we refer to
Crampin and Pirani \cite{crampin_pirani86}. Consider $L: (x,y) \in
TM \mapsto L(x,y)\in {\mathbb R}$ a first order regular Lagrangian
and let
\begin{equation}
\theta_L=d_JL=\frac{\partial L}{\partial y^i}dx^i, \
\omega_L=dd_JL=d\left(\frac{\partial L}{\partial y^i}\right)\wedge
dx^i \label{theta_omega} \end{equation} be the Cartan-Poincar\'e one
and two-forms. For these we have the following properties.

The vertical distribution of the tangent bundle is a Lagrangian
subbundle for the symplectic structure $\omega_L$. This is due to
the fact that $d_J\theta_L=d^2_JL=0$ and therefore
$\omega_L\left(JX, JY\right)=0,$ for all $X,Y$ vector fields on
$TM$.

For an arbitrary semispray $S$ on $TM$, which is a second order
vector field on the base manifold $M$, we have $i_S\theta_L={\mathbb
C}(L)$, where ${\mathbb C}=y^i\left({\partial }/{\partial
y^i}\right)$ is the Liouville vector field on $TM$. By taking the
exterior differential of the above formula we obtain
$i_S\omega_L={\mathcal L}_S\theta_L-d\left({\mathbb C}(L)\right).$
For the Euler-Lagrange vector field $S$ we have ${\mathcal
L}_S\theta_L=dL$. Therefore, the equation
$i_S\omega_L=-d\left(L-{\mathbb C}(L)\right)$ uniquely determines
the canonical semispray and implies that $S$ is a Cartan symmetry
for the Euler-Lagrange equations.

\section{Third order vector fields and Craig-Synge covectors}

Spaces with metric structures whose components depend on higher
order accelerations were introduced by Kawaguchi,
\cite{kawaguchi32}. The classic variational problem of a metric
structure of order two induces the Euler-Lagrange covector field,
which is of order four. From the corresponding vector field, which
lives on the tangent bundle of order three, it is difficult to
derive geometric objects that live on the tangent bundle of second
order and to develop a geometric theory of the given Lagrangian of
second order.

For a second order regular Lagrangian, using the Euler-Lagrange
covector field, Craig \cite{craig35} and Synge \cite{synge35}
derived another covector of order three along a curve. Miron
\cite{miron94} paid a special attention to this covector that was
called Craig-Synge covector. The corresponding vector field lives on
the tangent bundle of order two, it is a semispray and one can
develop from it the geometry of the given Lagrangian.

In this section we prove, that the Craig-Synge covector can be
independently derived using a variational principle of a regular
second order Lagrangian, for the case when only the vertical
components of the curve are varied. The corresponding vector field
is a semispray associated to the second order Lagrange space. In
this section we provide the equation that uniquely determines this
canonical semispray of a second order regular Lagrangian. This
equation generalizes the well-known equation for a first order
Lagrangian.

Consider the following variational problem for the regular
Lagrangian $L$. Let
\[ c_{\varepsilon}: t\in [0,1] \mapsto
c_{\varepsilon}\left(t\right)=\left(x^{i}\left(t\right),
\frac{dxi}{dt}\left(t\right) + \varepsilon V{}^{i}\left(t\right),
\frac{1}{2}\frac{d^2 x^i}{dt^2}\left(t\right)+
\varepsilon\frac{dV^{i}}{dt}\left(t\right)\right)\in T^{2}M,\] be a
variation of the curve $c\left(t\right)=c_{0}\left(t\right)$ in
$T^2M$, where $\varepsilon$ belongs to some small neighborhood of
$0\in\mathbb{R}$, and $V^{i}\left(t\right)$ are the components of a
vector field along curve $c$ such that
$V^{i}\left(0\right)=V^{i}\left(1\right)=0$. We look for necessary
conditions for the curve $c=c_{0}$ to be an extremal of the integral
\[
I\left(c_{\varepsilon}\right)=
\int_{0}^{1}L\left(x^{i}\left(t\right),
\frac{dx^i}{dt}\left(t\right) + \varepsilon V{}^{i}\left(t\right),
\frac{1}{2}\frac{d^2 x^i}{dt^2}\left(t\right) +
\varepsilon\frac{dV^{i}}{dt}\left(t\right)\right) dt.\] For this we
require that $c$ is a solution of the following equation
\begin{equation}
\left.\frac{d}{d\varepsilon}
\left(I\left(c_{\varepsilon}\right)\right)\right|_{\varepsilon=0}=0.
\label{eq:ddepsilon}\end{equation} Equation \eqref{eq:ddepsilon} can
be written as follow
\begin{multline*}
0=\left.\frac{d}{d\varepsilon}\left(I\left(c_{\varepsilon}\right)\right)\right|_{\varepsilon=0}
= \int_{0}^{1}\left(\frac{\partial L}{\partial
y^{\left(1\right)i}}V{}^{i} +
\frac{\partial L}{\partial y^{\left(2\right)i}} \frac{dV^{i}}{dt}\right)dt\\
=\int_{0}^{1}\left[\frac{\partial L}{\partial y^{\left(1\right)i}} -
\frac{d}{dt}\left(\frac{\partial L}{\partial
y^{\left(2\right)i}}\right)\right] V^{i}dt + \left.
\left(\frac{\partial L}{\partial
y^{\left(2\right)i}}V^{i}\right)\right|_{t=0}^{t=1}. \end{multline*}
Since $V$ is an arbitrary vector field, we have that equation
\eqref{eq:ddepsilon} holds true if and only if the following
Craig-Synge equations, \cite{craig35}, \cite{synge35} are satisfied:
\begin{equation}
\frac{\partial L}{\partial
y^{\left(1\right)i}}-\frac{d}{dt}\left(\frac{\partial L}{\partial
y^{\left(2\right)i}}\right)=0.\label{eq:craig-synge}
\end{equation}
For a regular Lagrangian $L$, equations (\ref{eq:craig-synge})
represent a system of third order differential equations, which can
be written as follows:
\begin{equation}
\frac{d^{3}x^{i}}{dt^{3}} + 3!
G^{i}\left(x,\frac{dx}{dt},\frac{1}{2!}\frac{d^{2}x}{dt^{2}}\right)=0.
\label{eq:k+1_system}
\end{equation} The functions $G^{i}$ are given
by \begin{equation}
3G^{j}=\frac{1}{2}g^{ji}\left[d_{T}\left(\frac{\partial L}{\partial
y^{\left(2\right)i}}\right)-\frac{\partial L}{\partial
y^{\left(1\right)i}}\right]\label{eq:canonical_gi}\end{equation} and
they are local coefficients of a semispray, as it has been shown by
Miron in \cite{miron94}. We refer to the vector field $S$ given by
expression (\ref{eq:sgi}), whose coefficients are given by
expression \eqref{eq:canonical_gi}, as to the semispray $S$ of the
regular Lagrangian $L$. In expression \eqref{eq:canonical_gi},
$d_{T}$ is the Tulczyjew operator, \cite{tulczyjew76}\[
d_{T}=y^{\left(1\right)i}\frac{\partial}{\partial
x^{i}}+2y^{\left(2\right)i}\frac{\partial}{\partial
y^{\left(1\right)i}}.\]

Next theorem gives the equation that uniquely determines the
semispray of a second order regular Lagrangian.

\begin{thm}
\label{thm:lstheta} For a regular second order Lagrangian $L$, the
semispray, whose coefficients are given by expression
\eqref{eq:canonical_gi}, is the only semispray $S$ on $T^{2}M$ that
satisfies the following equation\begin{equation}
\mathcal{L}_{S}\left(d_{J^{2}}L\right)=d_{J}L.\label{eq:lstheta}\end{equation}
\end{thm}
\begin{proof}
Consider $S$ a semispray on $T^{2}M$, given by expression \eqref{eq:sgi}.
The Lie derivative of the Cartan-Poincar\'e one-form $\theta_{L}^{2}=d_{J^{2}}L$
has the following expression
\begin{equation}
\mathcal{L}_{S}\theta_{L}^{2} = \left[d_{T}\left(\frac{\partial
L}{\partial y^{\left(2\right)i}}\right) -
6g_{ji}G^{j}\right]dx^{i}+\frac{\partial L}{\partial
y^{\left(2\right)i}}dy^{\left(1\right)i}.
\label{eq:lstheta_local}
\end{equation}
Using expression \eqref{eq:lstheta_local}, we obtain that for a
regular Lagrangian, equation \eqref{eq:lstheta} uniquely determines
the functions $G^{i}$, which are given by expression
\eqref{eq:canonical_gi}.
\end{proof}
\begin{cor}
For a regular second order Lagrangian $L$, its semispray is uniquely
determined by the following equation
\begin{equation}
i_{S}\omega_{L}^{2} = - d\left({\mathbb C}_{1}\left(L\right)\right)
+ \theta_{L}^{1}. \label{eq:isomega}
\end{equation}
\end{cor}
\begin{proof}
From Proposition \ref{pro:thetap} it follows that
$\mathcal{L}_{S}\theta_{L}^{2}=i_{S}\omega_{L}^{2}+d\left({\mathbb
C}_{1}\left(L\right)\right)$. Therefore, equations
\eqref{eq:lstheta} and \eqref{eq:isomega} are equivalent and each of
them, uniquely determine the semispray of the regular, second order
Lagrange space.
\end{proof}

For a first order regular Lagrangian $L(x,y)$, as it has been shown
by Sarlet et al. \cite{sarlet84}, the canonical semispray is
uniquely determined by the following equation
\begin{equation}
\mathcal{L}_{S}\left(d_{J}L\right)=dL.\label{eq:lsdjl}\end{equation}
For a first order Lagrangian the corresponding equation for
\eqref{eq:isomega} reads as follows
\begin{equation}
i_{S}\omega_L=-d\left({\mathbb
C}(L)-L\right).\label{eq:isomega1}\end{equation}

We remark here that in equations \eqref{eq:lsdjl} and
\eqref{eq:isomega1} the right hand-side is an exact form, while in
equations \eqref{eq:lstheta} and \eqref{eq:isomega} the right
hand-side is not an exact form any more since $\theta_{L}^{1}$ is
not even a closed one-form. This has the following consequences,
which are different from the case of a first order Lagrangian.

\begin{cor}
For a regular second order Lagrangian, its semispray satisfies:
\end{cor}
\begin{enumerate}
\item $\mathcal{L}_{S}\omega_{L}^{2}=\omega_{L}^{1}$, \medskip{}

\item $S\left({\mathbb C}_{1}\left(L\right)\right)={\mathbb C}_{2}\left(L\right)$.
\end{enumerate}
\begin{proof}
If we take the exterior differential of both sides of equation
\eqref{eq:lstheta} and use the fact that the Lie derivative commutes
with the exterior differential, we obtain
$\mathcal{L}_{S}\omega_{L}^{2}=\omega_{L}^{1}$. Using the skew
symmetry of $\omega_{L}^{2}$ and expression \eqref{eq:isomega}, we
obtain $S\left({\mathbb C}_{1}\left(L\right)\right)={\mathbb
C}_{2}\left(L\right)$.
\end{proof}
From the above results we can see that the semispray $S$ is not a
symmetry for the corresponding second order Lagrange space.

\section{Second order dynamical derivative}

In this section, we introduce dynamical covariant derivatives of
first and second order, induced by a pair $\left(S,N\right)$, where
$S$ is the associated semispray of a second order Lagrange space and
$N$ is a nonlinear connection. Using these dynamical derivatives and
the presymplectic structure we give conditions that uniquely
determine a nonlinear connection $N$. With respect to this nonlinear
connection, the second order dynamical covariant derivative of the
metric tensor vanishes and the horizontal distribution is a
Lagrangian subbundle for the presymplectic structure.

Let $S$ be the associated semispray of a second order Lagrange
space, which is uniquely determined by Theorem \ref{thm:lstheta} and
consider $N$ a nonlinear connection. For a $d$-vector field with
components $X^{i}$, we define its first and second order dynamical
derivatives as follows
\begin{equation}
\nabla X^{i} = S\left(X^{i}\right)+
\underset{\left(1\right)}{M_{j}^{i}\ }X^{j},\ \nabla^{2}X^{i} =
S^{2}\left(X^{i}\right)+2\underset{\left(1\right)}{M_{j}^{i}\
}S\left(X^{j}\right)+2\underset{\left(2\right)}{M_{j}^{i}\
}X^{j}.\label{eq:nablax}
\end{equation}

\begin{rem*}
Let $X=X^{i}\left(\partial/\partial x^{i}\right)$ be a vector field
on the base manifold $M$. Its complete lift, given by
\begin{equation} X^{c}=X^{i}\frac{\partial}{\partial x^{i}}+
S\left(X^{i}\right)\frac{\partial}{\partial
y^{\left(1\right)i}}+\frac{1}{2}S^{2}\left(X^{i}\right)\frac{\partial}{\partial
y^{\left(2\right)i}},\label{eq:xc1}\end{equation} can be expressed
as follows\begin{equation} X^{c}=X^{i}\frac{\delta}{\delta
x^{i}}+\nabla X^{i}\frac{\delta}{\delta
y^{\left(1\right)i}}+\frac{1}{2}\nabla^{2}X^{i}\frac{\partial}{\partial
y^{\left(2\right)i}}.\label{eq:xc2}\end{equation}

\end{rem*}
One can extend these first and second order dynamical derivatives to
arbitrary $d$-tensor fields. In this paper we are interested in the
first and second derivatives of the metric tensor $g_{ij}$. These
derivatives are given by
\begin{equation}
\nabla g_{ij}=S\left(g_{ij}\right)-
\underset{\left(1\right)}{N_{i}^{k}\
}g_{kj}-\underset{\left(1\right)}{N_{j}^{k}\
}g_{ik},\label{eq:nabla1g}
\end{equation}
\begin{eqnarray}
\nabla^{2}g_{ij} & = & S^{2}\left(g_{ij}\right)
- 2\underset{\left(1\right)}{N_{i}^{k}\ } S\left(g_{kj}\right)
- 2\underset{\left(1\right)}{N_{j}^{k}\ }S\left(g_{ik}\right)\nonumber \\
&  & -2\underset{\left(2\right)}{N_{i}^{k}\ }g_{kj}-
2\underset{\left(2\right)}{N_{j}^{k}\
}g_{ik}+2g_{mk}\underset{\left(1\right)}{N_{j}^{k}\
}\underset{\left(1\right)}{N_{i}^{m}\
}\label{eq:nabla2g}
\end{eqnarray}
Consider $g=g_{ij}dx^{i}\otimes dx^{j}$ the metric tensor of a
second order regular Lagrangian. One can extend this tensor to the
following metric structure on $T^{2}M$
\begin{eqnarray}
g^{c} & = & g_{ij}\left(dx^{i}\otimes dy^{\left(2\right)j} +
dy^{\left(1\right)i}\otimes dy^{\left(1\right)j}+
dy^{\left(2\right)i}\otimes dx^{j}\right)+\nonumber \\
&  & S\left(g_{ij}\right)\left(dx^{i}\otimes
dy^{\left(1\right)j}+dy^{\left(1\right)i}\otimes
dx^{j}\right)+\frac{1}{2}S^{2}\left(g_{ij}\right)dx^{i}\otimes
dx^{j}.\label{eq:gc1}\end{eqnarray}

With respect to the adapted basis of the nonlinear connection $N$,
the metric structure $g^{c}$ has the following expression
\begin{eqnarray}
g^{c} & = g{}_{ij}\left(dx^{i}\otimes\delta y^{\left(2\right)j} +
\delta y^{\left(1\right)i}\otimes\delta y^{\left(1\right)j}+
\delta y^{\left(2\right)i}\otimes dx^{j}\right)+\nonumber \\
& \nabla g_{ij}\left(dx^{i}\otimes\delta y^{\left(1\right)j}+\delta
y^{\left(1\right)i}\otimes
dx^{j}\right)+\frac{1}{2}\nabla^{2}g_{ij}dx^{i}\otimes
dx^{j}.\label{eq:gc2}
\end{eqnarray}

\begin{thm}
\label{thm:canonicaln} Let $S$ be the associated semispray of a
second order regular Lagrangian. There exists a unique nonlinear
connection $N$ on $T^{2}M$ such that
\end{thm}
\begin{enumerate}
\item $\omega_{L}^{2}\left(hX,v_{1}Y\right)=2g^{c}\left(hX,v_{1}Y\right),$
$\forall X,Y\in\chi\left(T^{2}M\right).$ \medskip{}

\item $\nabla^{2}g=0$, $d_{h}\theta_{L}^{2}=0.$
\end{enumerate}
\begin{proof}
We prove that the first condition of the theorem uniquely determines
the first coefficients $\underset{\left(1\right)}{N_{i}^{k}\ }$,
while the second condition uniquely determines the second
coefficients $\underset{\left(2\right)}{N_{i}^{k}\ }$ of a nonlinear
connection. We use the metric tensor $g_{ij}$ to raise and lower
indices, therefore we denote \[ \underset{\left(1\right)}{N_{ij}\
}:= g_{ik}\underset{\left(1\right)}{N_{j}^{k}\
},\underset{\left(2\right)}{N_{ij}\
}:=g_{ik}\underset{\left(2\right)}{N_{j}^{k}\ }.\]

The presymplectic structure $\omega_{L}^{2}$ can be expressed, with
respect to the adapted basis of a nonlinear connection $N$, as
follows
\begin{equation} \omega_{L}^{2}=\frac{\delta}{\delta
x^{i}}\left(\frac{\partial L}{\partial
y^{\left(2\right)j}}\right)dx^{j}\wedge dx^{i}+\frac{\delta}{\delta
y^{\left(1\right)i}}\left(\frac{\partial L}{\partial
y^{\left(2\right)j}}\right)\delta y^{\left(1\right)j}\wedge
dx^{i}+2g_{ij}\delta y^{\left(2\right)j}\wedge
dx^{i}.\label{eq:omega2ad}
\end{equation}
Condition 1 of the theorem can be written as follows
\begin{equation} \omega_{L}^{2}\left(\frac{\delta}{\delta
x^{i}},\frac{\delta}{\delta
y^{\left(1\right)j}}\right)=\frac{\delta}{\delta
y^{\left(1\right)i}}\left(\frac{\partial L}{\partial
y^{\left(2\right)j}}\right)=2g^{c}\left(\frac{\delta}{\delta
x^{i}},\frac{\delta}{\delta
y^{\left(1\right)j}}\right)=2g_{ij|},\label{eq:2gijbar}\end{equation}
which is equivalent to
\begin{equation}
2\underset{\left(1\right)}{N_{ij}\
}+4\underset{\left(1\right)}{N_{ji}\
}=2S\left(g_{ij}\right)+\frac{\partial^{2}L}{\partial
y^{\left(1\right)i}\partial
y^{\left(2\right)j}}.\label{eq:2n4n}\end{equation} Since $\nabla
g_{ij}=g_{ij|}$ are the components of a second rank symmetric
$d$-tensor field, from expression \eqref{eq:2gijbar} we obtain the
symmetry in $i$ and $j$ of the following expression \[
\frac{\delta}{\delta y^{\left(1\right)i}}\left(\frac{\partial
L}{\partial y^{\left(2\right)j}}\right),\] which is equivalent to
\begin{equation} 2\underset{\left(1\right)}{N_{ij}\
}-2\underset{\left(1\right)}{N_{ji}\ }=\frac{\partial^{2}L}{\partial
y^{\left(2\right)i}\partial
y^{\left(1\right)j}}-\frac{\partial^{2}L}{\partial
y^{\left(2\right)j}\partial
y^{\left(1\right)i}}.\label{eq:n1minusn1}\end{equation} From
expressions \eqref{eq:2n4n} and \eqref{eq:n1minusn1} we obtain that
$\underset{\left(1\right)}{N_{ij}\ }$ is uniquely determined and it
has the following expression \begin{equation}
\underset{\left(1\right)}{N_{ij}\
}=\frac{1}{3}S\left(g_{ij}\right)+\frac{1}{3}\frac{\partial^{2}L}{\partial
y^{\left(2\right)i}\partial
y^{\left(1\right)j}}-\frac{1}{6}\frac{\partial^{2}L}{\partial
y^{\left(2\right)j}\partial
y^{\left(1\right)i}}.\label{eq:n1ij}\end{equation} Using expression
\eqref{eq:nabla2g}, we obtain that first part of condition 2 of the
theorem, $\nabla^{2}g_{ij}=0$, uniquely determines the symmetric
part of $\underset{\left(2\right)}{N_{ij}\ }$.\begin{equation}
2\underset{\left(2\right)}{N_{ij}\
}+2\underset{\left(2\right)}{N_{ji}\
}=S^{2}\left(g_{ij}\right)-2S\left(g_{ik}\right)\underset{\left(1\right)}{N_{j}^{k}\
}-2S\left(g_{kj}\right)\underset{\left(1\right)}{N_{i}^{k}\
}+2g_{mk}\underset{\left(1\right)}{N_{j}^{k}\
}\underset{\left(1\right)}{N_{i}^{m}\
}.\label{eq:symn2ij}\end{equation} Second part of condition 2 of the
theorem can be written as follows
$\omega_{L}^{2}\left(hX,hY\right)=0,\forall
X,Y\in\chi\left(T^{2}M\right).$ From expression \eqref{eq:omega2ad}
this condition implies the symmetry in $i$ and $j$ of the following
expression \[ \frac{\delta}{\delta x^{i}}\left(\frac{\partial
L}{\partial y^{\left(2\right)j}}\right),\] which determines the
skewsymmetric part of $\underset{\left(2\right)}{N_{ij}\ }.$
\begin{equation}
2\underset{\left(2\right)}{N_{ij}\ } -
2\underset{\left(2\right)}{N_{ji}\ }=\frac{\partial^{2}L}{\partial
y^{\left(2\right)i}\partial x^{j}}-\frac{\partial^{2}L}{\partial
y^{\left(2\right)j}\partial
x^{i}}-\underset{\left(1\right)}{N_{j}^{k}\
}\frac{\partial^{2}L}{\partial y^{\left(2\right)i}\partial
y^{\left(1\right)k}}+\underset{\left(1\right)}{N_{i}^{k}\
}\frac{\partial^{2}L}{\partial y^{\left(2\right)j}\partial
y^{\left(1\right)k}}.\label{eq:n2minusn2}\end{equation} Using
expressions \eqref{eq:symn2ij} and \eqref{eq:n2minusn2}, we obtain
that the coefficients $\underset{\left(2\right)}{N_{ij}\ }$ are
uniquely determined.
\end{proof}
\begin{cor}
With respect to the nonlinear connection determined by Theorem
\ref{thm:canonicaln}, the presymplectic structure
$\omega_{L}^{2}=d\theta_{L}^{2}$ has the following form
\begin{equation} \omega_{L}^{2}= d\theta_{L}^{2}=2g_{ij|}\delta
y^{\left(1\right)j}\wedge dx^{i}+2g_{ij}\delta
y^{\left(2\right)j}\wedge
dx^{i}.\label{eq:adaptedomega}\end{equation}

If $h$, $v_{1}$ and $v_{2}$ are the projectors that correspond to
the nonlinear connection determined by Theorem \ref{thm:canonicaln},
then the Poincar\'e-Cartan one -form $\theta_{L}^{2}$ satisfies the
following equations\begin{equation}
d_{h}\theta_{L}^{2}=d_{v_{1}}\theta_{L}^{2} =
d_{v_{2}}\theta_{L}^{2}=0.\label{eq:dhtheta2}\end{equation}
\end{cor}
\begin{proof}
First condition and second part of the second condition of Theorem
\ref{thm:canonicaln} are equivalent with the fact that the presymplectic
structure $\omega_{L}^{2}$ has with respect to adapted basis of the
nonlinear connection the form given by expression \eqref{eq:adaptedomega}.
All equations \eqref{eq:dhtheta2} can be obtained directly from expression
\eqref{eq:adaptedomega} if we compose it with the projectors $h$,
$v_{1}$ and $v_{2}$ respectively.
\end{proof}
Theorem \ref{thm:canonicaln} gives conditions that uniquely
determine a nonlinear connection for a second order Lagrange space.
These conditions can be restated as follows, there is a unique
nonlinear connection such that $\nabla^{2}g=0$ and the presymplectic
structure $\omega_{L}^{2}$ is given by expression
\eqref{eq:adaptedomega}. However, the expressions we obtained for
the coefficients $\underset{\left(1\right)}{N_{ij}\ }$ and
$\underset{\left(2\right)}{N_{ij}\ }$in the proof of the above
theorem are not easy to use. Next proposition gives a simpler form
for the first coefficients of the considered nonlinear connection.

\begin{prop} \label{prop:n1ij}
Consider $S$ the associated semispray of a second order regular
Lagrangian $L$, whose coefficients $G^{i}$ are given by expression
\eqref{eq:canonical_gi}. The following functions are the first
coefficients of the nonlinear connection determined by Theorem
\ref{thm:canonicaln}. \begin{equation}
\underset{\left(1\right)}{N_{j}^{i}\ } =\frac{\partial
G^{i}}{\partial y^{\left(2\right)j}},\ \label{eq:n1}\end{equation}
\end{prop}
\begin{proof}
We have to show that the functions
$\underset{\left(1\right)}{N_{j}^{i}\ }$ given by expression
\eqref{eq:n1} satisfy expression \eqref{eq:n1ij}, which uniquely
determine the first coefficients of the nonlinear connection. The
associated semispray $S$ of a second order regular Lagrangian has
the local coefficients given by expression \eqref{eq:canonical_gi}
and it is uniquely determined by the following condition
\begin{equation} S\left(\frac{\partial L}{\partial
y^{\left(2\right)i}}\right)=\frac{\partial L}{\partial
y^{\left(1\right)i}}.\label{eq:sl2}\end{equation} By direct
calculation we have the following formula\begin{equation}
\left[S,\frac{\partial}{\partial
y^{\left(2\right)i}}\right]=-2\frac{\partial}{\partial
y^{\left(1\right)i}}+3\frac{\partial G^{j}}{\partial
y^{\left(2\right)i}}\frac{\partial}{\partial
y^{\left(2\right)j}}.\label{eq:sdy2}\end{equation} If we apply both
sides of expression \eqref{eq:sdy2} to $\partial L/\partial
y^{\left(2\right)k}$ and use expression \eqref{eq:sl2} we obtain
\begin{equation} g_{jk}\frac{\partial G^{j}}{\partial
y^{\left(2\right)i}}=\frac{1}{3}S\left(g_{ik}\right)+\frac{1}{3}\frac{\partial^{2}L}{\partial
y^{\left(2\right)k}\partial
y^{\left(1\right)i}}-\frac{1}{6}\frac{\partial^{2}L}{\partial
y^{\left(2\right)i}\partial
y^{\left(1\right)k}}.\label{eq:pgpy2}\end{equation} From expressions
\eqref{eq:n1ij} and \eqref{eq:pgpy2} we obtain that first
coefficients $\underset{\left(1\right)}{N_{j}^{i}\ }$ of the
nonlinear connection are given by expression \eqref{eq:n1}.
\end{proof}
According to Proposition \ref{prop:n1ij} the coefficients
$\underset{(1)}{N^i_j}$ of the nonlinear connection determined by
Theorem \ref{thm:canonicaln} are the same for the nonlinear
connections that were considered by Catz \cite{catz74}, Dodson and
Radivoiovici \cite{dodson82}, de Leon and Rodrigues \cite{deLeon85},
Byrnes \cite{byrnes96}, Miron \cite{miron97} and Bucataru
\cite{bucataru98}. However, the other coefficients
$\underset{(2)}{N^i_j}$ of the nonlinear connection we obtained in
Theorem \ref{thm:canonicaln} are different and with respect to them
the presymplectic form has a very simple form given by expression
\eqref{eq:adaptedomega}.

\section{Examples}

In this section we consider examples of second order regular
Lagrangians, which were previously used by Miron and Atanasiu
\cite{mironat96} to extend semi-Riemannian, Finslerian and
Lagrangian structures from the base manifold to the total space of
second order tangent bundles. For such second order regular
Lagrangians we compute the associated semispray and the nonlinear
connection as they have been defined in the previous sections. If
the Lagrangian function depends on position effectively, this
nonlinear connection is different from the prolongation of a
nonlinear connection introduced by Catz \cite{catz74} and used by
Miron and Atanasiu \cite{mironat96}. It is also different from the
nonlinear connection used by Dodson and Radivoiovici \cite{dodson82}
and de León and Vasquez \cite{deLeon85} to study the geometry of the
tangent bundle of order two. Each of these nonlinear connections
give information regarding the Lagrangian function as we will see in
this section. However, for the semi-Riemannian case, with respect to
the nonlinear connection introduced in the previous section, the
first and second order covariant derivatives of the metric tensor
vanish and the presymplectic structure has a very simple form with
respect to it.

Let $g=g_{ij}\left(x\right)dx^{i}\otimes dx^{j}$ be a
semi-Riemannian metric on the base manifold $M$ and denote by
$\gamma_{jk}^{i}\left(x\right)$ the Christoffel symbols of the
Riemannian metric. Consider the first order Lagrangian $L_1:TM
\longrightarrow \mathbb{R}$ given by
\begin{equation}
L_1\left(x,y^{(1)i}\right)=g_{ij}(x)y^{(1)i}y^{(1)j}=\left\|y^{(1)}\right\|^2.
\label{eq:l1} \end{equation} We recall that the geodesics of the
semi-Riemannian metric $g$, which are extremal curves for the first
order Lagrangian $L_1$, are solutions of the following system of
equations:
\begin{equation}
\nabla\left(\frac{dx^i}{dt}\right)=0. \label{eq:rgeod}
\end{equation}
Using the fact that $\nabla g=0$, for a curve $c(t)=(x^i(t),
y^{(1)i}(t))$ on $TM$ we have that
\begin{equation}
\frac{d}{dt}\left(L_1\left(x(t), y^{(1)i}(t)\right)\right)=
2g_{ij}(x(t))\frac{dx^i}{dt}\nabla y^{(1)i}, \label{eq:conservL1}
\end{equation} from which we obtain that $L_1$ is conserved along
the geodesics of the semi-Riemannian metric $g$. The right hand side
of expression (\ref{eq:conservL1}) is the first-order deformation
Lagrangian ${\mathcal L}_1$, considered by Casciaro and Fracaviglia
in \cite{casciaro97}, expression (3.3).

By direct calculation it follows that
\begin{equation}
z^{\left(2\right)i}=y^{\left(2\right)i}+
\frac{1}{2}\gamma_{jk}^{i}\left(x\right)
y^{\left(1\right)j}y^{\left(1\right)k} =\frac{1}{2}\nabla y^{(1)i}
\label{eq:z2i}\end{equation} are the components of a $d$-vector
field on $T^{2}M$, which represents the covariant expression of the
second order acceleration. The function
$L_2:T^{2}M\longrightarrow\mathbb{R}$, given by
\begin{equation}
L_2\left(x,y^{\left(1\right)},y^{\left(2\right)}\right)=
g_{ij}\left(x\right)z^{\left(2\right)i}z^{\left(2\right)j} =
\frac{1}{4}\left\|\nabla y^{(1)i}\right\|^2 \label{eq:riemann_L}
\end{equation}
is well defined on the second order tangent space, it
is a second order regular Lagrangian and represents the square of
the magnitude of the second order acceleration. The following
formulas can be obtained by a straightforward calculation.
\begin{equation} \frac{\partial
L_2}{\partial y^{\left(2\right)i}}=2g_{ij}z^{\left(2\right)j},
\frac{\partial L_2}{\partial y^{\left(1\right)i}}=
2g_{kj}z^{\left(2\right)j} \gamma_{ip}^{k}y^{\left(1\right)p},
\frac{\partial^{2}L_2}{\partial y^{\left(1\right)i}\partial
y^{\left(2\right)j}} = 2g_{kj} \gamma_{ip}^{k}y^{\left(1\right)p}.
\label{eq:dldy}\end{equation} The associated semispray of the second
order regular Lagrangian given by expression \eqref{eq:riemann_L} is
uniquely determined by Theorem \ref{thm:lstheta}. Its coefficients
are given by \begin{eqnarray} 3G^{i} & = &
d_{T}\left(z^{\left(2\right)i}\right) +
\gamma_{jk}^{i}z^{\left(2\right)j}y^{\left(1\right)k} \nonumber \\
& = & \frac{1}{2}\left(\frac{\partial\gamma{}_{jk}^{i}}{\partial
x^{m}}+ \gamma_{pj}^{i} \gamma_{km}^{p}\right) y^{\left(1\right)j}
y^{\left(1\right)k} y^{\left(1\right)m} +3\gamma_{jk}^{i}
y^{\left(1\right)j}y^{\left(2\right)k}.
\label{eq:riemann_3gi}\end{eqnarray} The corresponding nonlinear
connection of the second order Lagrangian $L_2$ is uniquely
determined by Theorem \ref{thm:canonicaln}. First set of
coefficients of this nonlinear connection are given by expression
\eqref{eq:n1}, which in our case are given by
\begin{eqnarray} \underset{\left(1\right)}{N_{j}^{i}\
}\left(x,y^{\left(1\right)}\right) =\frac{\partial G^{i}}{\partial
y^{\left(2\right)j}}= &
\gamma_{jk}^{i}\left(x\right)y^{\left(1\right)k}\label{eq:riemann_n1}\end{eqnarray}
With respect to the nonlinear connection studied in
\cite{bucataru98} and using the second order covariant derivative,
the Craig-Synge equations, given by expression
\eqref{eq:k+1_system}, of the second order Lagrangian $L_2$ can be
expressed in a very simple form, which generalizes the geodesic
equation \eqref{eq:rgeod}
\begin{equation}
\nabla^2\left(\frac{dx^i}{dt}\right)=0. \label{eq:nebla2xi}
\end{equation}
With respect to the nonlinear connection introduced by Miron and
Atanasiu in \cite{mironat96}, for a curve $c(t)=\left(x^i(t),
dx^i/dt \right)$ on $TM$ we have the following formula
\begin{equation}
\frac{d^2}{dt^2}\left(L_1\left(x, \frac{dx}{dt}\right)\right)=
8L_2\left(x, \frac{dx}{dt}, \frac{1}{2}\frac{d^2 x}{dt^2}\right) +
2g_{ij}(x(t))\frac{dx^i}{dt}\nabla^2 \left( \frac{dx^i}{dt}\right),
\label{eq:conservL2}
\end{equation}
which is equivalent to
\begin{equation}
\frac{d^2}{dt^2}\left\|\frac{dx}{dt} \right\|^2= 2\left\|\nabla
\left( \frac{dx^i}{dt} \right) \right\|^2 +
2g_{ij}(x(t))\frac{dx^i}{dt}\nabla^2 \left( \frac{dx^i}{dt} \right).
\label{eq:norm}
\end{equation}

With respect to the nonlinear connection introduced in the previous
section the presymplectic structure and the complete lift have a
very simple form. From expression \eqref{eq:adaptedomega} we obtain
that the canonical presymplectic structure of the second order
Lagrange space can be expressed as follows\begin{equation}
\omega_{L_2}^{2}=2g_{ij}{\delta y}^{\left(2\right)j}\wedge
dx^{i}.\label{eq:riemann_omega2L}\end{equation} From expression
\eqref{eq:gc2}, the complete lift $g^{c}$ of the Riemannian metric
tensor can be expressed as follows:
\begin{equation}
g^{c}=g_{ij}\left(dx^{i}\otimes \delta y^{\left(2\right)j}+ \delta
y^{\left(1\right)i}\otimes \delta y^{\left(1\right)j}+ \delta
y^{\left(2\right)i}\otimes
dx^{j}\right).\label{eq:gc_riemann}\end{equation}

In this section we have seen that the covariant derivatives induced
by the considered semispray and the nonlinear connection give simple
covariant expressions for the geometric structures associated to the
second order Lagrangian $L_2$.

\begin{acknowledgement*}
This work was supported by grants PN II IDEI 398 (I.B.) and by CEEX
252 (R.M.) of the Romanian Ministry of Education.
\end{acknowledgement*}

\end{document}